\documentstyle[11pt]{article}    
\title{Singular connection and Riemann theta function}  %

\author{Weiping Li}   
\date{ }
\begin{document}           

\newtheorem{lm}{Lemma}[section]
\newtheorem{thm}[lm]{Theorem}
\newtheorem{pro}[lm]{Proposition}
\newtheorem{cor}[lm]{Corollary}
\newtheorem{df}[lm]{Definition}
\newcommand{\sect}[1]{\section{#1} \setcounter{equation}{0}}
\renewcommand{\theequation}{\arabic{section}.\arabic{equation}}
\def\eot{{\hfill \vrule height7pt width10pt}\\ \vspace{3 ex} }
\maketitle                 

\begin{abstract}
We prove the Chern-Weil formula for $SU(n+1)$-singular connections over
the complement of an embedded oriented surface in smooth four manifolds.
The expression of the representation of a number as a sum of nonvanishing
squares is given in terms of the representations of a number as a sum of
squares. Using the number theory result, we study the irreducible
$SU(n+1)$-representations of the fundamental group of the complement 
of an embedded oriented surface in smooth four manifolds.


\end{abstract}

\sect{Introduction}

We study $SU(n+1)$-singular connections over $X \setminus \Sigma$ in
this paper, where $X$ is a smooth closed oriented 4-manifold and $\Sigma$
is a closed embedded surface. In \cite{sw}, S. Wang first started to understand
the topological information from sigular connections. Later, Kronheimer and 
Mrowka \cite{km} studied the Donaldson invariants under the change of 
$SU(2)$-singular connections. The paper \cite{km} turns out to be a crucial 
step for analyzing the structure of the Donaldson invariants and for
recent development in the Seiberg-Witten theory.

In \S 2, we describe the $SU(n+1)$-singular connection space over 
$X \setminus \Sigma$. The $SU(n+1)$-singular Chern-Weil formula is given in
Proposition~\ref{cw}. In order to study the irreducible 
$SU(n+1)$-representations, we need to study the sum of nonvanishing squares.
This is one of well-known number theory problems. Jacobi initially studied the
representation $r_k(n)$ of a number $n$ as a sum of $k$-squares via 
Riemann theta function as a generating function. For a topological reason, we 
would like to understand the representation ${\cal R}_k(n)$ of a 
number $n$ as a sum of nonvanishing $k$-squares. In general it is difficult to
calculate both $r_k(n)$ and ${\cal R}_k(n)$. We prove a nice relation between
$r_k(n)$ and ${\cal R}_k(n)$ in Proposition~\ref{rR}. Up to the author's 
knowledge, it has not known before how to express the number ${\cal R}_k(n)$
(see \cite{gro}). Similar relation for the representation of a number by a 
quadratic form is also obtained.

In the last section, we use those number theory criterions to study 
$SU(n+1)$-singular flat connections. This gives an interesting interaction
between ${\cal R}_n(N)$ in the number theory and the irreducible
$SU(n+1)$-representations of $\pi_1( X \setminus \Sigma)$ in topology 
(see Proposition~\ref{inter}).

\sect{$SU(n+1)$ singular connection and Chern-Weil formula}

\noindent {\bf (i) Singular connections over $X \setminus \Sigma $}
\vspace{0.2in}

Let $X$ be a smooth closed oriented 4-manifold and let
$\Sigma$ be a closed embedded
surface. We will assume that both $X$ and $\Sigma $ are connected, and $\Sigma$
to be orientable or oriented for simplifying our discussions. 
Denote $tN$ be a closed tubular neighborhood of $\Sigma$. 
Identify $tN$ diffeomorphically to the unit disk bundle of the normal bundle. 
Let $Y$ be the boundary of $tN$, 
which has the circle bundle structure over $\Sigma $ by this diffeomorphism. 
Let $i \eta$ be a connection 1-form for the circle bundle, 
so $\eta $ is an $S^1$-invariant 
1-form on $Y$ which coincides with the 1-form $d \theta $ on each 
circle fiber. Using $(r, \theta )$ polar coordinates in some 
local trivialization of the disk bundle, 
we have that $dr \wedge d \theta$ is the correct orientation 
for the normal plane. 
By radial projection, $\eta $ can be extended to $tN \setminus  \Sigma$.
\vspace{0.1in}

We will  work on the  structure group $SU(n+1)$ for $n \geq 1$. 
So a connection $A$ on $X \setminus \Sigma $ which 
is represented on each normal plane to 
$\Sigma $ by a connection matrix looks like
\begin{equation}\label{21}
i \left( \begin{array}{cccc}
\alpha_0 & & & \\
	 & \alpha_1 &  & \\
	 &       & \ddots & \\
	 &    &  & \alpha_n  \end{array}
\right) d \theta + (\mbox{lower terms)}, \ \ \ 
\sum^n_{i=0} \alpha_i \equiv 0 \pmod 1 .  \end{equation}
The size of the connection matrix is $o(r^{-1})$, so $A$ is singular along 
the surface $\Sigma $.
\vspace{0.1in}

For every $SU(n+1)$-singular connection, one can associate with holonomy
as in \cite{ss, sw}. 
Let $P \to  X \setminus \Sigma$ be a vector bundle with 
structure group $SU(n+1)$ for $n \geq 1$. 
To define the holonomy around $\Sigma $
of a connection $A$ on $P$, for any point 
$\sigma \in \Sigma$ and real number $0 < r <1$, let 
$S^1_{\sigma}(r)$ be a circle 
with center $\sigma $ and radius $r$ on the normal plane of $tN$ over $
\sigma$. An element $h(A; \sigma, r) \in SU(n+1)$ can be obtained by parallel 
transport of a frame of $P$ along $S^1_{\sigma }(r)$ via the connection $A$. 
Although $h(A; \sigma, r)$ depends on the choice of a frame, 
its conjugacy class
$[h(A; \sigma, r)]$ in $SU(n+1)$ does not (c.f. \cite{ss, sw}). 
If for all $\sigma \in \Sigma$, the limit 
$h_A = \lim_{r \to 0^+}[h(A; \sigma, r)]$ is independent of 
$\sigma $ and $tN$, we call it the holonomy of $A$ along $\Sigma$.
\vspace{0.1in}

The holonomy of this connection on the positively oriented small circles of 
radius $r$ is approximately
\begin{equation} \exp 2 \pi i \left( \begin{array}{cccc}
\alpha_0 & & & \\
	 & \alpha_1 &  & \\
		  &       & \ddots & \\
			   &    &  & \alpha_n  \end{array}
			   \right)
, \ \ \ \sum^n_{i=0} \alpha_i \equiv 0 \pmod 1 .\end{equation} 
Since  only the conjugacy class of the  holonomy has any invariant meaning, we
may suppose that $\alpha_i$ lies in the interval $[0, 1)$,
therefore the matrices (2.2) modulo the permutation group $S_{n+1}$ 
run through each conjugacy  class just once.
When $\alpha_i = 0$ for all $0 \leq i \leq n$, the holonomy is trivial,
and if the phrase ``lower terms'' makes sense, we have ordinary connections on
$X$. Also when $\alpha_i \in \{ 0, \frac{1}{n}, \cdots, \frac{n-1}{n}\}$
for all $i$, 
the holonomy is in the center of $SU(n+1)$ (n-th root of unity), 
the associated $SU(n+1)/{Z_n}$-bundle has trivial holonomy; 
and with this twist
we can consider these as $PSU(n+1)$-connections on $X$.
\vspace{0.1in}

Conjugacy classes in $SU(n+1)$ can be characterized by parameters
$\alpha_i$ with 
$$\alpha  = (\alpha_i )_{0 \leq i \leq n} \in [0, 1)^{n+1}
/(\alpha_0 + \alpha_1 + \cdots + \alpha_n \equiv 0 \pmod 1).$$
Note that any permutation of $(\alpha_i )$ gives the same conjugacy class.
Hence we can stay on the region for conjugacy classes by making $\alpha_i$
satisfying the following
\begin{equation}
1 > \alpha_0 \geq \alpha_1 \geq \cdots \geq \alpha_n \geq 0,
\ \ \mbox{and} \ \  \alpha_0 + \alpha_1 + \cdots 
+ \alpha_n \in \{0, 1, \cdots, n\}. \end{equation} 
The region of conjugacy classes of $SU(n+1)$ can be identified with
the quotient space $\{z_i \in S^1, \prod_{i=0}^n z_i = 1\}/S_{n+1}$
of the maximal torus of $SU(n+1)$ under the Weyl group action.
When $n = 1$ ($SU(2)$ case), $1 > \alpha_0 \geq \alpha_1 \geq 0, \alpha_0 
+ \alpha_1 = 0, \alpha_0 + \alpha_1 = 1$. The equation $\alpha_0
+ \alpha_1 = 0$ always gives the identity conjugacy class. So  
$1 > \alpha_0 \geq \alpha_1 \geq 0$ and $\alpha_0 + \alpha_1 = 1$ describe
the conjugacy classes of $SU(2)$ as in \cite{km} with
$\alpha_0^{'} = \alpha_0 + 1/2, \alpha_1^{'} = \alpha_1 + 1/2$.
For $\alpha_i = 0 \geq \alpha_{i+1} \geq \cdots \alpha_n \geq 0$, the
conjugacy classes can be viewed as conjugacy classes in $SU(i)$
or $U(i)$ in $SU(n+1)$.

The matrix-valued 1-form given on $X \setminus  \Sigma$ by the expression 
\begin{equation}  i \left( \begin{array}{cccc}
\alpha_0 & &  & \\
    & \alpha_1 & & \\
    &      & \ddots  & \\
    &      &       & \alpha_n  
    \end{array}  \right) \eta , \end{equation} 
has the asymptotic behavior of $(2.1)$, but is only defined locally. To make an
$SU(n+1)$ connection on $X \setminus  \Sigma$ which has this form near $\Sigma$
, start with  $SU(n+1)$ bundle $\overline{P}$ over $X$ and choose a $C^{\infty}$
decomposition of $\overline{P}$ on $N$ as 
\[\overline{P}|_N = \overline{L_0}  \oplus \overline{L_1}  \oplus \cdots 
\oplus \overline{L_n}, \]
compatible with the hermitian metric. 
Since we will work on the modeled connection $A^{\alpha}$, the
decomposition of $\overline{P}|_N$ gives a natural model.
Although  $\overline{P}|_N$ is trivial, 
but $\overline{L_i}$ may not be. We define
topological invariants in this situation:
\begin{equation} \left \{  \begin{array}{ll}
k& =c_2(\overline{P})[X] \\
l_i  & = -c_1(\overline{L_i})[\Sigma ], \ \ \ \sum^n_{i=0}l_i = 0.
\end{array}  \right. \end{equation} 
Choose any smooth $SU(n+1)$  connection $\overline{A}^0$ on 
$\overline{P}$ which
respects to the decomposition over $N$, so  we have 
\begin{equation} \overline{A}^0|_N = \left( \begin{array}{cccc}
b_0& &  & \\
 & b_1 & & \\
  &      & \ddots  & \\
      &      &       &b_n
      \end{array}  \right), \ \ \ \sum^n_{i=0}b_i = 0, \end{equation} 
where  $b_i$ is a smooth connection in $\overline{L_i}$. 
Let the model
connection $A^{\alpha}$ on $P =  \overline{P}|_{X \setminus  \Sigma}$ be
the following:
\begin{equation}
A^{\alpha} = \overline{A}^0 + i \beta(r)  \left( \begin{array}{cccc}
\alpha_0 & &  & \\
    & \alpha_1 & & \\
	&      & \ddots  & \\
	    &      &       & \alpha_n  
		\end{array}  \right) \eta , \end{equation}
where $\beta$ is a smooth cutoff function equal to 1 in $[0, \frac{3}{8}]$
and equal to 0 for $r \geq \frac{1}{2}$. In terms of trivialization compatible 
with the decomposition, the second term in $A^{\alpha}$ is an element  of 
$\Omega^1_{N \setminus \Sigma }(Ad P)$ which can be extended to all of 
$X \setminus \Sigma$. The curvature $F(A^{\alpha })$ extends to a smooth 2-form
with values in $Ad \overline{P}$ on the whole $X$, since $i d \eta $ is 
smooth on $tN$,  $i \eta$ is the pullback to $tN$ of the 
curvature form of the circle bundle $Y$. It can be
thought as a smooth 2-form on the surface $\Sigma $.
\vspace{0.1in}

The connection $A^{\alpha}$ in (2.7) defines a connection on $X\setminus \Sigma$
.  The holonomy $h_{A^{\alpha}}$  around small 
linking circles is asymptotically equal to $(2.2)$. 
We now define an affine space of connections modeled on $A^{\alpha}$ 
by choosing $p > 2$ and denoting
\[{\cal A}^{\alpha, p}_1 = \{ A^{\alpha } + 
a \ \ | \ \|a\|_{L^P(X \setminus \Sigma)} + 
\|\nabla_{A^{\alpha}} a\|_{L^p{(X \setminus \Sigma)}} < \infty  \}.\]
Similarly a gauge group
\[{\cal G}^{\alpha, p}_2 
= \{ g \in Aut P \ \ | \ \ \|g\|_{L^p(X \setminus \Sigma)}
+ \|\nabla_{A^{\alpha}}g\|_{L^p(X \setminus \Sigma)} +
\|\nabla_{A^{\alpha}}^2g\|_{L^p(X \setminus \Sigma)}< \infty  \}.\]
The $L^p$ space is defined by using the measure 
inherited from any smooth measure on $X$.

\begin{pro}  (i) The space ${\cal A}^{\alpha, p}_1$ and 
${\cal G}^{\alpha, p}_2$ are independent of the 
choices of $\overline{A}^0$ and the connection 1-form $\eta$.

(ii) The space ${\cal G}^{\alpha, p}_2$ is a Banach Lie group 
which acts smoothly on ${\cal A}^{\alpha, p}$ and is independent 
of $\alpha $. The stabilizer of $A$
is $Z_n$  or $H$ ($Z_n \subset H \subset SU(n+1)$)
according  as $A$ is irreducible or reducible respectively.
\end{pro}
Proof: The proof is same as in the Proposition 2.4 
in \cite{km} (see also Chapter 3 \cite{sw}). 
\eot

\noindent {\bf (ii) The Chern-Weil formula for ${\cal A}^{\alpha ,p}_1$}
\vspace{0.1in}

By the same token of the
Proposition 3.7 in \cite{km}, we have that the equivalence class of the norm $
L^p_{k, A^{\alpha}}$ is independent of $\alpha =  (\alpha_0, \alpha_1, \cdots,
\alpha_n)$, so the gauge group 
${\cal G}^{\alpha, p}_{k+1} = {\cal G}^p_{k+1}$ is independent
of $\alpha $ as a parameter in the definition of the 
model connection $A^{\alpha}$.
So the space ${\cal A}^{\alpha, p}_1 = A^{\alpha } + L^p_1(\Omega (X 
\setminus \Sigma ), Ad P))$, where $L^p_1(\Omega (X 
\setminus \Sigma ), Ad P))$
is a Banach space which is independent of $\alpha $, and $a \in
L^p_1(\Omega (X \setminus \Sigma ), Ad P))$, the diagonal component $D(a)$ of 
$a$ is in
$L^p_1$, and $(a - D(a))$ is in $L^p_1$, $r^{-1} (a - D(a))$ is in $L^p$.
\vspace{0.1in}

\begin{pro} \label{cw} 
For all $A \in {\cal A}^{\alpha, p}$, the following 
Chern-Weil formula holds.
\begin{equation} \label{cwf}
\frac{1}{8 \pi^2} \int_{X \setminus \Sigma  } tr (F_A \wedge F_A) = 
k +  \sum^n_{i=0} \alpha_i l_i - \frac{1}{2} (\sum^n_{i=0} \alpha_i^2)  
\Sigma \cdot \Sigma  . \end{equation}
\end{pro} 
Proof: We begin by proving the formula for the model connection  $A^{\alpha}$
in the simplest  case. Let $A^{\alpha}$ be globally reducible. So $\overline{E}
=  \overline{L_0} \oplus \overline{L_1} \cdots \oplus \overline{L_n}$ globally, that $\overline{b_i}$  is  a smooth connection on $\overline{L_i}$ and that $A^{\alpha}$ is reducible as
\[ A^{\alpha } =  \left( \begin{array}{cccc}
b_0 & & &  \\
    & b_1 &  & \\
    & \ddots & &  \\
    &  & & b_n  
    \end{array} \right) \ \ \ \ \mbox{with} \ \ \ \ 
    b_i = \overline{b_i} + i  \alpha_i  
    \beta(r) \eta .  \]
The cutoff function $\beta (r)$ is defined in (2.7), and $i \eta $ is a 
connection 1-form on the normal circle bundle to $\Sigma $.
The closed 2-form $d(i \beta(r) \eta)$ can be extended smoothly across $\Sigma$,
due to $\beta(r) =1$ near $\Sigma $. So it is the pullback  from $\Sigma $
of the  curvature  form $F(i \eta)$. Since the second cohomology of the 
neighborhood is 1-dimensional, so the closed 2-form  
$d(i \beta(r) \eta)= F(i \eta)$ represents  the Poincar\'{e} 
dual of $\Sigma $, denote  by  $P.D(\Sigma ) = d(i \beta(r) \eta)$.
Hence we have the degree of the normal bundle  
\[ \int_{\Sigma } - \frac{1}{2 \pi i} F(i \eta ) = 
\int_{\Sigma } - \frac{1}{2 \pi i} d(i \beta(r) \eta) =  < P.D(\Sigma ), 
\Sigma >  =  \Sigma \cdot \Sigma . \]
In de Rham cohomology, we have
 \begin{eqnarray}
 - \frac{1}{2 \pi i} d(b_i) &  = & - \frac{1}{2 \pi i} d(\overline{b_i}) +  
\alpha_i (- \frac{1}{2 \pi i}d(i \beta(r) \eta))  \nonumber \\
 &  = &  c_1(\overline{L_i}) + \alpha_i  P.D(\Sigma ) \end{eqnarray} 
So we have the following two identities:
\begin{equation} \label{32}
\frac{1}{8 \pi^2}  tr(F_{A^{\alpha}} \wedge F_{A^{\alpha}}) = 
\frac{1}{8 \pi^2}  tr (dA^{\alpha} \wedge dA^{\alpha}) =  
\frac{1}{8 \pi^2} tr (\sum^n_{i=0}
db_i \wedge db_i) , \end{equation}
\begin{equation} \label{33}
< (- \frac{1}{2 \pi i} db_i)  \wedge (- \frac{1}{2 \pi i} db_i),  X > =  -  <
  \frac{1}{4 \pi^2}  db_i \wedge db_i  , X > . \end{equation}
Therefore by (\ref{33}) and (\ref{32}),
\begin{eqnarray}
 - \frac{1}{4 \pi^2}  < d b_i \wedge d b_i, [X] > & = & 
- < (c_1(\overline{L}_i) + \alpha_i P.D(\Sigma )) \wedge
(c_1(\overline{L}_i) + \alpha_i P.D(\Sigma )), [X] > \nonumber \\
& = & - < (c_1(\overline{L}_i))^2, [X] >  - 2 \alpha_i < (c_1(\overline{L}_i)
\wedge P.D(\Sigma )), [X] > \nonumber \\
& & - (\alpha_i )^2 < P.D(\Sigma ) \wedge P.D(\Sigma ),
[X] > \nonumber \\
& = & - l_i^2 - 2 \alpha_i < c_1(\overline{L}_i),  \Sigma  > - 
(\alpha_i )^2 < P.D(\Sigma ), X \cap \Sigma > \nonumber \\
& = & - l_i^2 + 2 \alpha_i l_i  - \alpha_i^2 \Sigma \cdot \Sigma . 
\end{eqnarray} 
Observe that on the Lie algebra ${\bf su}(n)$ of skew adjoint matrices
$tr (M^2) = - |M|^2$.
Hence by (\ref{32}) and (\ref{33}) the Chern-Weil formula 
for the modeled connection $A^{\alpha }$ is 
\begin{equation} \label{35}
\frac{1}{8 \pi^2} \int_{X \setminus \Sigma} tr F_{A^{\alpha }} \wedge 
F_{A^{\alpha }} = \frac{1}{2}\sum_{i=0}^n (- l_i^2 + 2 \alpha_i l_i  
- \alpha_i^2 \Sigma \cdot \Sigma). \end{equation}
Note that $c_2(\overline{E}) = \sum_{i < j}c_1(\overline{L}_i) \cdot c_1(
\overline{L}_j) = \sum_{i < j} l_i l_j$. Also from
$c_1(\overline{E}) = \sum_{i=0}^n l_i = 0$, we have 
\[ 0 = (\sum_{i=0}^n l_i)^2 = \sum_{i=0}^n l_i^2 + 2  \sum_{i < j} l_i l_j, \]
i.e.  $ c_2 (\overline{E}) = - \frac{1}{2} \sum_{i=0}^n l_i^2.$ 
By (\ref{35}) we have
\begin{equation} \label{36}
\frac{1}{8 \pi^2} \int_{X \setminus \Sigma} tr F_{A^{\alpha }} \wedge 
F_{A^{\alpha }} = k  + \sum^n_{i=0} \alpha_i \cdot l_i - \frac{1}{2}
(\sum_{i=0}^n\alpha_i^2) (\Sigma \cdot \Sigma ).\end{equation}
Although this calculation is global, it has an interpretation locally on $tN$.
Let $Y_{\varepsilon } \subset tN$ be the 3-manifold circle bundle over $\Sigma $
given by radius $r = \varepsilon $. The Chern-Simons integral is given by
the following:
\[ cs_{\varepsilon }(A^{\alpha}) = \frac{1}{8 \pi^2} \int_{Y_{\varepsilon }} 
tr(dA^{\alpha} \wedge A^{\alpha} + \frac{2}{3}A^{\alpha} \wedge A^{\alpha}
\wedge A^{\alpha}). \]
The integral $cs_{\varepsilon }(A^{\alpha})$ depends only on the 
homotopy class of the trivialization of 
the bundle on $Y_{\varepsilon}$ with respect to which
the connection matrix $A^{\alpha}$
is computed (see \cite{dk}). Since there is a distinguished trivialization on 
$Y_{\varepsilon }$ which extends to $tN$, we have the 
Chern-Simons $cs_{\varepsilon }$ defined as a real number. 
Let $X_{\varepsilon}$ be
the complement of $\varepsilon $-neighborhood of $\Sigma $ with boundary
$Y_{\varepsilon}$. Thus
\[\frac{1}{8 \pi^2} \int_{X_{\varepsilon}} tr F_A \wedge F_A = k + 
cs_{\varepsilon }(A^{\alpha }) . \]
By (\ref{36}) from the reducible connection, we have
\[ \lim_{\varepsilon \to 0} cs_{\varepsilon }(A^{\alpha }) = 
\sum^n_{i=0} \alpha_i \cdot l_i - \frac{1}{2}
(\sum_{i=0}^n\alpha_i^2 ) \Sigma \cdot \Sigma .\]
So the Chern-Weil formula holds whenever $A$ is a connection which 
is smooth and reducible near to $\Sigma $ by applying the above local
statement. Since such connections are dense in
${\cal A}^{\alpha, p}_1$ and the curvature integral is a continuous function of 
$A$ in the $L^p_{1,A^{\alpha}}$-topology, the result follows.
\eot

\noindent {\bf Remarks:} (1) When $\alpha_i = 0$ for all $i$, 
Proposition~\ref{cw} is the usual Chern-Weil formula. 

(2) When $n=1$, we have the $SU(2)$-situation. 
Proposition~\ref{cw} for $\alpha_i = \alpha^{'} +1/2 (i=0, 1)$
and $l_0 + l_1 = 0$ case is 
\[ \frac{1}{8 \pi^2} \int_{X \setminus \Sigma} tr F_A \wedge F_A = 
k + 2 (\alpha_0^{'}) l_0 - (\alpha_0^{'})^2 (\Sigma \cdot \Sigma ) , \]
which is the Proposition 5.7 in \cite{km}. 
So our formula  
extends their formula to the $SU(n+1)$ group.
\vspace{0.1in}

\begin{cor} Let $a$ be the restriction of 
$A \in {\cal A}^{\alpha , p}_1$ on the
boundary of $X \setminus \Sigma $. Then the Chern-Simons invariant takes 
the value
\begin{equation}  cs(a) \equiv \sum_{i=0}^n \alpha_i l_i - 
\frac{1}{2}(\sum_{i=0}^n \alpha_i^2
)  \Sigma \cdot \Sigma  \ \ \ \pmod 1. \end{equation} 
\end{cor} \eot

The proof follows directly from the proof of Proposition~\ref{cw}. 
The Chern-Weil formula gives the charge for singular $SU(n+1)$-connections
over $X \setminus \Sigma $. We study the maximum and minimum of the 
charge over the conjugacy holonomy region.

\begin{cor} For $\Sigma \cdot \Sigma \neq 0$, the charge 
takes its maximum and minimum by comparing $k$ with the following values
\[k + \frac{\sum_{i=0}^n l_i^2}{2\Sigma \cdot \Sigma }
- \frac{j^2 \Sigma \cdot \Sigma }{2(n+1)},  \ \ \ \ j = 1, 2, \cdots, n;\]
\[ k + \frac{\sum_{i=0}^{j-1} l_i^2}{2\Sigma \cdot \Sigma }
- \frac{m^2 \Sigma \cdot \Sigma }{2j} + \frac{s_j}{j}(m - \frac{s_j}
{2\Sigma \cdot \Sigma }), \ \ \ j = 2, 3, \cdots, n; m= 1, \cdots, j-1;\]
where $s_j = \sum_{i=0}^{j-1} l_i$ (not necessary zero).
\end{cor} 
Proof: First of all, we find out the extreme values inside the region. By the
method of 
Lagrange multiplier, 
\[f(\alpha_0, \cdots, \alpha_n) = k + \sum_{i=0}^n \alpha_i l_i 
- \frac{1}{2}(\sum_{i=0}^n \alpha_i^2) \Sigma \cdot \Sigma ,\]
with constraints $\sum_{i=0}^n \alpha_i = j (j=1, \cdots, n)$,
so we have the critical point 
\[ \alpha_i =  \frac{l_i}{\Sigma \cdot \Sigma } + \frac{j}{n+1}, 
\ \ \ \ \ 0 \leq i \leq n,\] and its corresponding
charge is 
\[k+ \frac{\sum_{i=0}^n l_i^2}{\Sigma \cdot \Sigma } - 
\frac{j^2 \Sigma \cdot \Sigma }{2 (n+1)}, \ \ \ \ j =1, \cdots, n.\]
For $j=0$, all $\alpha_i$ are zero, so the charge is $k$.

By the method of Lagrange multiplier to locate all critical points of 
$f(\alpha_0, \cdots, \alpha_n)$ on the boundary $\alpha_j = \cdots = \alpha_n = 0 (j=2, \cdots, n)$ of the region $(2.3)$, we have the function
\[ f(\alpha_0, \cdots, \alpha_{j-1}; \lambda_m) = 
k + \sum_{i=0}^{j-1} \alpha_i l_i 
- \frac{1}{2}(\sum_{i=0}^{j-1} \alpha_i^2) \Sigma \cdot \Sigma 
- \lambda_m (\sum_{i=0}^{j-1} \alpha_i - m), \]
with constraints $\sum_{i=0}^{j-1} \alpha_i = m$ for $m =1, \cdots, j-1$.
The critical point is
\[ \alpha_i = \frac{l_i}{\Sigma \cdot \Sigma } + \frac{m}{j} -
\frac{\sum_{i=0}^{j-1}l_i}{j \Sigma \cdot \Sigma },  \ \ \ \
i = 0, 1, \cdots, j-1 ,\]
and its corresponding charge is, by a straightforward calculation, 
\[k + \frac{\sum_{i=0}^{j-1} l_i^2}{2\Sigma \cdot \Sigma }
- \frac{m^2 \Sigma \cdot \Sigma }{2j} + \frac{s_j}{j}(m - \frac{s_j}
{2\Sigma \cdot \Sigma }), \]
where $s_j = \sum_{i=0}^{j-1} l_i$. 
So the result follows by comparing these extreme values. \eot

\sect{Riemann theta function}

In this section, the needed number theory criteria are shown.
The Riemann theta function enters our picture naturally from the 
representation of a number as a sum of $k$-squares, or by a quadratic form.
It is one of well-known classic problems in number theory. Jacobi initially 
studied this problem by using Riemann theta function as a generating function;
in particular, Siegel  generalized vastly to several complex variables
(see \cite{mu}). 
At the moment we are only interested in the
relation between topology and number theory (see \S 4). Further 
investigation along this line will be discussed elsewhere.
\vspace{0.1in}

\noindent {\bf (i) The representation of a number as a sum  of nonvanishing
squares.}
\vspace{0.1in}

Let $N$ be an integer, $N \geq 1$, with the representation 
\[ N = N_1^2 + N_2^2 + \cdots + N_n^2 , \]
where the $(N_i)_{1 \leq i \leq n}$'s are integers including zero. Let $r_n
(N)$ denote the number of representations of $N$ as the sum of $n$ squares.
\begin{enumerate}
\item For $n =2$, Jacobi derived an identity from the generating function
\[\theta_3(0,z) = \sum^{\infty }_{n = - \infty} q^{n^2}, \ \ \ 
q = e^{\pi i z}, \ \ \mbox{with} \ \ Im z > 0.  \]
Jacobi identity  is 
\[\{\theta_3(0,z) \}^2 = 1 + 4 \sum_{n =1}^{\infty } \frac{(-1)^{n-1} q^{2n-1}}
{1 - q^{2n-1}}, \]
which gives the result that
\[r_2(N) = 4 \sum_{d odd, d|N} (-1)^{\frac{1}{2}(d-1)} = 4 \{d_1(N) - d_3(N)\}.\]
where $d_1(N)$ and $d_3(N)$ are the numbers of the divisors of $N$ of the form 
$4m +1$ and $4m+3$ respectively.

\item For $n=3$, Legendre proved that a number $N$ is the sum of three squares
if and only if $N \neq 4^a (8b + 7), a \geq 0 , b \geq 0$. For all
$N$, $r_3(4^a N) = r_3(N)$. The function $r_3(N)$
has been evaluated by Dirichlet as a finite sum involving symbols of quadratic
reciprocity. We give the following formula for $r_3(N)$ (see \cite{gro}):
\[ r_3(N) = \frac{G_N}{\pi} \sqrt{N} L(1, \chi ), \]
where \[G_N = \left\{ \begin{array}{ll}
0 & N \equiv 0, 4, 7 \pmod 8 \\
16 & N \equiv 3 \pmod 8 \\
24 & N \equiv 1, 2, 5, 6 \pmod 8
\end{array} \right. \]
and $L(s, \chi ) = \sum_{m=1}^{\infty} \chi (m) m^{-s}$ with 
$\chi (m) = (-4N/m)$, $(r/N)$ the Jacobi symbol.

\item For $n \geq 4$, Lagrange  (1770) proved that every positive integer can be
represented as the sum of four squares, hence also as the sum of five or more
squares. In particular for $n=4$,
\[r_4(N) = 8 \sum_{d|N, 4\not| d} d , \]
where the summation is over those positive divisors of $N$, which are not 
divisible by 4.
\end{enumerate}

For our purpose, we need to get the representation $N$ as the sum of
nonvanishing 
integer squares. Let ${\cal R}_n(N)$ denote 
the number of representations of $N$ as 
the sum of $n$ nonvanishing squares. For example, $r_1(3) = r_2(3) = 0$, 
$r_3(3) = 8$, $r_4(3) = 32$, but 
${\cal R}_4(3) = 0$ (see (\ref{rR1}). 
Note that $r_1(N) = {\cal R}_1(N)$ for all $N$.
The following proposition gives the relation between $r_n(N)$ 
and ${\cal R}_n(N)$.
\begin{pro} \label{rR}
If $N$ is an integer ($N \geq 1$), then
\begin{equation} \label{rR1}
{\cal R}_n(N) = \sum_{i=1}^n (-1)^{n-i}\left( \begin{array}{c} n \\ i
\end{array}  \right)  r_i(N) . \end{equation}
\end{pro} 
Proof:  If $q =  e^{\pi i z}$ with  $Im z > 0$, then  by definition we have
\[\{ \theta_3(0,z)\}^n = (\sum_{l = - \infty }^{\infty } q^{l^2})^n =
\sum_{N=0}^{\infty } r_n(N)q^N, \ \ \ r_n(0) = 1.\]
The generating function for ${\cal R}_n(N)$ is $\theta_3(0,z) -1$, so 
\[\{\theta_3(0,z) -1\}^n  = (\sum_{l \neq 0}q^{l^2})^n = \sum_{N=1}^{\infty}
{\cal R}_n(N) q^N. \]
On the other hand the binomial formula gives
\[\{\theta_3(0,z) -1\}^n  = 
\sum_{i=0}^n (-1)^{n-i} \left( \begin{array}{c} n \\ i
\end{array}  \right)\{ \theta_3(0,z)\}^i . \]
For the constant coefficient, it corresponds to $i=0$, i.e. $(1-1)^n = 0$.
Therefore by comparing the coefficients of $q^N (N \geq 1)$, we have 
the desired relation. \eot

Although we know that for each $n \geq 5$ all but a finite set of integers
are sums of exactly $n$ nonvanishing squares (see \cite{gro} Chapter 6), 
Proposition~\ref{rR} gives the precise relation among the
numbers of representations of $N$ as sums of squares and sums of
nonvanishing squares.
\vspace{0.1in}

\noindent {\bf (ii) The representation of a number by a quadratic form}
\vspace{0.1in}

Let $(a_{pq})$ be a real, symmetric, $n \times n$ matrix, 
and let the associated quadratic form $Q(x) = \sum_{p,q = 1}^n a_{pq}x_p x_q$ 
be positive definite. 
It is well-known  that the multiple  series  
\[ \sum_{i_1, \cdots, i_n = - \infty}^{\infty} e^{\pi i z Q(i_1, \cdots, i_n)}\]
converges absolutely and uniformly in every compact set in the upper half-plane
$Im z > 0$. The theta function associated  to $Q$ is defined to be 
\[ \theta(z,Q) =  \sum_{i_1, \cdots, i_n = - \infty}^{\infty} 
e^{\pi i z Q(i_1, \cdots, i_n)} . \]
In case $a_{pq} = \delta_{pq}$ is the identity matrix, then the $\theta(z, Id)$ 
reduces to $\{\theta_3(0,z)\}^n$. In our application later, we have the matrix
even, i.e. $a_{pp}$ are even. Then the definition of $\theta(z,Q)$  yields
\[ \theta (z +1, Q) = \theta (z, Q) . \]
In the next section we will consider the particular even matrix:
\begin{equation}
(a_{pq}) = \left( \begin{array}{cccc}
2  & 1 & \cdots & 1 \\
1 & 2 & \cdots & 1 \\
\vdots & \vdots & \ddots  & \vdots \\
1 & 1 & \cdots &  2 
\end{array} \right)_{n \times n}  . \end{equation} 
Its  determinant is $n+1$.
\begin{thm}
Let $(a_{pq})$ be a symmetric, $n \times n$ matrix of integers, where $a_{pp}$
are all even for $p = 1, 2, \cdots, n$, and the associated quadratic form 
$Q(x)$ be positive definite with determinant $D$. 
Let $Q^{-1}$ be the inverse form of $Q$. Then we have  
\[\theta (z+1, Q) = \theta (z, Q), \ \ \ 
\theta (- \frac{1}{z}, Q) = (\sqrt{\frac{z}{i}})^n D^{- \frac{1}{2}} \theta (z,
Q^{-1}), \]
for all complex  $z$ with $Im z > 0$.
\end{thm}

From the above relations, one can derive the formula for $\theta (\frac{az +b}
{cz +d}, Q)$, with $a, b, c, d$ are integers and $ad - bc = 1$, since the 
modular group is generated by the two transformations  $A: z \to z +1$, and
$B: z \to - \frac{1}{z}$ (see \cite{mu}).
\vspace{0.1in}

Let $r_Q(N)$ (or ${\cal R}_Q(N)$) denote the 
number of (or all nonzero) solutions 
$x_1, \cdots, x_n$, with $x_i$ integral for every $i$, 
such that $1 \leq i \leq n$
of the equation
\[ \sum_{p, q = 1}^n a_{pq} x_p x_q = 2 N . \]

Let $(a_{pq})_i$ be the $(n-1) \times (n-1)$ matrix by deleting i-th row and
i-th column of the matrix $(a_{pq})$. Denote the corresponding quadratic form
be $Q_i$. Clearly $Q_i$ is an even, symmetric, positive definite form.
Similarly $Q_{i_1i_2}$ is the quadratic form with $x_{i_1} = x_{i_2} = 0$, etc.
The following lemma gives the relation
among $r_Q(N)$, $r_{Q_{i_1 \cdots i_j}}(N)
(j=1, 2, \cdots, n-1)$ and ${\cal R}_Q(N)$.

\begin{pro} \label{033}
For the even quadratic form $Q$, we have the relation
\[{\cal R}_Q(N) = r_Q(N) - \sum_{i=1}^n r_{Q_i}(N) + \]
\begin{equation} \label{330}
\sum_{1\leq i_1 < i_2 \leq n}r_{Q_{i_1i_2}}(N) - \cdots
+ (-1)^{n-1} \sum_{1\leq i_1 < \cdots < i_{n-1} \leq n} r_{Q_{i_1\cdots i_{n-1}}}(N) . \end{equation}
\end{pro}
Proof: Since $Q(x)$ is an even form and $\theta(z, Q)$ is holomorphic for $Imz
> 0$, so an expansion of $\theta (z, Q)$ in power of $e^{2\pi i z}$ is given 
by 
\[\theta (z, Q) = 1 + \sum_{N =1}^{\infty } r_Q(N) e^{2 \pi i N z}, \ \ \ 
Im z > 0 .\]
There is an another way to write the expansion of $\theta (z,Q)$ as
\begin{eqnarray}  
\theta (z, Q) & = & 1 + \sum^{\infty }_{N =1} {\cal R}_Q(N) e^{2 \pi i N z}
+ \sum^n_{i =1} (  \sum_{x_i =0} e^{\pi i z Q_i(x)})
- \sum^n_{i =1}(  \sum_{x_i= x_j =0} e^{\pi i z Q_i(x)}) + \cdots \nonumber \\ 
& = & 1 + \sum^{\infty }_{N =1} {\cal R}_Q(N) e^{2 \pi i N z} +\\
& & \sum^n_{i =1} (\sum^{\infty }_{N =1}r_{Q_i}(N) e^{2 \pi i N z})
- \sum^n_{i =1} (\sum^{\infty }_{N =1}r_{Q_{i_1i_2}}(N)e^{2 \pi i N z}) +
\cdots .   \nonumber
\end{eqnarray}
Hence the relation follows by comparing the coefficients of $e^{2 \pi i N z}$.
\eot

In particular, if we take the matrix $Q = 2 Id$, then 
\[r_Q(N) = r_n(N), \ \ \ \mbox{and} \ \ \ {\cal R}_Q(N) = {\cal R}_n(N) . \]
So the above lemma gives Proposition~\ref{rR}. It is clear from our discussion
that there are general relations between the number of solutions and the number
of nonvanishing solutions via the recursive formula in the theta function.

\sect{Unitary representation of $\pi_1 (X \setminus \Sigma )$}

In this section, we will use previous results to derive the nontrivial 
$SU(n+1)$-representations of $\pi_1 (X \setminus \Sigma )$. 

\begin{lm} \label{41}
Let $A$ be a flat $\alpha$-twisted $SU(n+1)$ connection. Then the holonomy 
parameter $(\alpha_i )_{0 \leq i \leq n}$, the instanton number  $k$ and 
monopole numbers  $l_i$ are related by
\begin{eqnarray}
l_i & = & \alpha_i (\Sigma \cdot \Sigma ), \ \ \ \mbox{for} \ \ \
0 \leq i \leq n  , \\
k& = & - \frac{\sum_{i=0}^n l_i^2}{2 \Sigma \cdot \Sigma } .
\end{eqnarray}
If $\Sigma \cdot \Sigma = 0$, then $k =0$ and $l_i =0$ for all $i$.
\end{lm} \label{4.1}
Proof: The flat $\alpha $-twisted connection
$A$ is one of the model connection corresponding to some integers
$k, l_i$. Since the bundle is flat, we have 
\[ w = diag ( c_1(\overline{L}_i) + \alpha_i 
P.D (\Sigma ))_{0 \leq i \leq n} \ \ \ ,\]
where the 2-form $w$ is in the proof of Proposition 2.2. Thus for each $i$,
$c_1(\overline{L}_i) + \alpha_i P.D (\Sigma ) = 0$; the equality $l_i 
=\alpha_i (\Sigma \cdot \Sigma )$ follows from integrating over $\Sigma $. If 
$\Sigma \cdot \Sigma = 0$, $l_i =0$ for all $i$ as well.
\vspace{0.1in}

On the other hand, if $A$ is flat, the Chern-Weil formula gives
\begin{eqnarray}
0 & = & \frac{1}{8 \pi^2} \int_{X \setminus \Sigma }tr F_A \wedge F_A \nonumber
\\
& = & k + \sum_{i=0}^n \alpha_i l_i - \frac{1}{2} (\sum_{i=0}^n \alpha_i^2 )
\Sigma \cdot \Sigma  \nonumber \\
& = & k + \frac{\sum_{i=0}^n l_i^2}{2 \Sigma \cdot \Sigma } \nonumber
\end{eqnarray}
If $\Sigma \cdot \Sigma = 0$, we have $k =0$ from the second equality. 
We obtain the formula (4.2). \eot

\noindent {\bf Remark:} Note that for $SU(n+1)$-flat bundles 
$|l_i | < | \Sigma \cdot \Sigma |$ due to $0 \leq \alpha_i <1$ from (4.1).
All $\alpha_i$'s are rational.
Using $\sum_{i=0}^n l _i = 0, l _0 = - \sum_{i=1}^n l_i$
 then we have 
\begin{equation} \label{qua}
k = - \frac{2 \sum_{i=1}^n l_i^2 + \sum_{i\neq j}^n l_i l_j}{2 
 \Sigma \cdot \Sigma }.
\end{equation} 

\begin{pro} \label{inter}
For a simply connected $X$ and an embedded oriented surface 
$\Sigma$ with $\Sigma \cdot \Sigma \neq 0$, 
\begin{enumerate}
\item If $\sum_{i,j = 1}^n l_i l_j = 0$, $\Sigma \cdot \Sigma$ is not a divisor
of any $N ( < n(\Sigma \cdot \Sigma)^2)$ with ${\cal R}_n(N) \neq 0$;
\item  In general $\Sigma \cdot \Sigma$ is not a divisor of any $N (<
\frac{n(n+1)}{2} (\Sigma \cdot \Sigma)^2)$ with ${\cal R}_Q(N)
\neq 0$ for $Q$ as (3.2); 
\end{enumerate}
then $\pi_1(X \setminus \Sigma )$ has no irreducible 
representation in $SU(n+1)$.
\end{pro}
Proof: Suppose there were an irreducible representation 
$\rho : \pi_1 (X \setminus
\Sigma ) \to SU(n+1)$ $(n \geq 1)$. The image of $\rho$ does not contain
in any proper subgroup of $SU(n+1)$. Denote $A$ be the corresponding flat 
connection on $X \setminus \Sigma $. By Seifert-Ven Kempf theorem, we have
\[ \begin{array}{ccl}
\pi_1 (Y_{\varepsilon }) & \to & \pi_1 (X \setminus \Sigma) \\
\downarrow &  & \downarrow \\
\pi_1(N_{\varepsilon }) & \to & \pi_1(X) = \{1\} .
\end{array}  \]
So the holonomy on $\pi_1(X \setminus \Sigma )$ is same as on 
$\pi_1(Y_{\varepsilon })$. The space
$Y_{\varepsilon }$ is the $S^1$-bundle over $\Sigma $, the homotopy exact  sequence of the fibration $S^1 \to Y_{\varepsilon } \to \Sigma $ yields
\[ \{1 \} \to \pi_1(S^1) \to \pi_1 (Y_{\varepsilon }) \to 
\pi_1(\Sigma ) \to \{1\} . \]
In other words, $\pi_1 (Y_{\varepsilon })$ is a central  extension of $
\pi_1(\Sigma )$. Let $\gamma $ be a generator of $\pi_1(S^1)$.

Since the conjugacy class $[\gamma ]$ generates $\pi_1(X \setminus \Sigma)
$ and $\rho $ is an irreducible representation, 
so the holonomy of $\rho_A$ is not in $Z_n$ and other proper subgroups.

Therefore $\rho_A : \pi_1 (X \setminus \Sigma) \to SU(n+1)/Z_n = PSU(n+1)$. 
\[ \pi_1(Y_{\varepsilon }) = \{ a_i, b_i, \gamma | 
\prod_i [a_i, b_i] = \gamma^m
, [\gamma^m, a_i] =1, [\gamma^m, b_i] =1 \}, \]
where $m = | \Sigma \cdot \Sigma |$ the absolute value of $\Sigma \cdot \Sigma$.
So the representation $\rho_A(\gamma )^m = \prod_i [\rho_A(a_i), \rho_A(b_i)]$,
and $[\rho_A(\gamma )^m, \rho_A(a_i)] = 1, [\rho_A(\gamma )^m, \rho_A(b_i)]
= 1$. This derives that the matrix $\rho_A(\gamma )$ must be a diagonal matrix: 
\[ \rho_A(\gamma ) = exp 2 \pi i \left( \begin{array}{cccc}
\alpha_0 & & & \\
& \alpha_1 & & \\
&   & \ddots & \\
& & & \alpha_n 
\end{array} \right), \]
for $\alpha_i$'s in the domain (2.3).
Now we can take the flat connection $A$ as 
our model connection corresponding to $\alpha = (\alpha_0, \cdots, \alpha_n)$
and instanton number $k$, monopole numbers $l_i$. By Lemma~\ref{41}
and (\ref{qua}), we have  
\begin{eqnarray}
l_i & = & \alpha_i \Sigma \cdot \Sigma   \nonumber \\
k & = & - \frac{2 \sum_{i=1}^n l_i^2 + \sum_{i \neq j}^n l_i l_j }{2 
\Sigma \cdot \Sigma } 
\end{eqnarray} 
If $\sum_{i,j= 1}^n l_i l_j = 0$, $k = - \frac{\sum_{i=1}^n l_i^2}{
\Sigma \cdot \Sigma }$ (an integer), then 
for $(1)$ (same argument for $(2)$) the resulting  
number on the right hand (4.4) is not an integer by the very definition of 
${\cal R}_n(N)$
(${\cal R}_Q(N)$) in \S 3 with $N = \sum_{i=1}^nl_i^2$
($N=  \sum_{i=1}^n l_i^2 + \frac{1}{2} \sum_{i,j= 1}^n l_i l_j$). 
Note that $l_i = \alpha_i \Sigma \cdot \Sigma 
< \Sigma \cdot \Sigma$, thus the range for $N = \sum_{i=1}^n l_i^2$ is
$< n (\Sigma \cdot \Sigma)^2$ ($< \frac{n(n+1)}{2} (\Sigma \cdot \Sigma)^2$). 

If $\Sigma \cdot \Sigma = \pm 1$, then any number $N < n (\pm 1)^2$ 
($N < \frac{n(n+1)}{2} (\pm 1)^2$) has
${\cal R}_n(N) = 0$ (${\cal R}_Q(N) = 0$). 
So at least one of $l_i = 0$ for $N = \sum_{i=1}^nl_i^2$ ($N = 
\sum_{i=1}^n l_i^2 + \frac{1}{2} \sum_{i,j= 1}^n l_i l_j$),
i.e. the corresponding $\alpha_i = 0$. Hence the induced image of $\rho_A$ is
a proper subgroup of $SU(n+1)$. So $\pi_1(X \setminus \Sigma)$ has no 
irreducible $SU(n+1)$-representations.
\eot

\noindent{\bf Remarks:}
\begin{enumerate}
\item We need to use the definition of ${\cal R}_n(N)$ to cover the
case of $\Sigma \cdot \Sigma = \pm 1$. For $n=1$ Proposition~\ref{inter}
is the Corollary 5.8 in \cite{km}.
\item The condition $\sum_{i,j= 1}^n l_i l_j = 0$ is different from $
\sum_{i,j= 0}^n l_i l_j = 0$. The later one with $\sum_{i=0}^n l_i =0$
will imply all $l_i =0$;
we have all $l_i \neq 0$, otherwise it will reduce to $SU(m)$ or $U(m) 
(m < n+1)$. 
So we may take $n+1$ as minimum number of $l_i \neq 0$.
\item One can (inductively) apply Proposition~\ref{inter} to ${\cal R}_{k}(N)$
for non 
representations of $\pi_1(X \setminus \Sigma )$ in an rank $k$ subgroup of
$SU(n+1)$.
\end{enumerate}

\vspace{0.5in}

Department of Mathematics
 
Oklahoma State University 

Stillwater, OK 74078. 

Email address: wli@littlewood.math.okstate.edu.

Fax: 405 - 744 - 8275 
\end{document}